# Reviewing Gödel's and Rosser's meta-reasoning of undecidability

Bhupinder Singh Anand

I review the classical conclusions drawn from Gödel's meta-reasoning establishing an undecidable proposition GUS in standard PA. I argue that, for any given set of numerical values of its free variables, every recursive arithmetical relation can be expressed formally in PA by different, but formally equivalent, propositions. I argue that this asymmetry yields alternative Representation and Self-reference meta-Lemmas. I argue that Gödel's meta-reasoning can thus be expressed avoiding any appeal to the truth of propositions in the standard interpretation IA of PA. I argue that this now establishes GUS as decidable, and PA as omega-inconsistent. I argue further that Rosser's extension of Gödel's meta-reasoning involves an invalid deduction.

## 1. Introduction

**1.1 The main argument**

I argue that, for any given set of natural number values of its free variables, every recursive arithmetical relation[1], say $q(x, y)$, can be expressed ([Me64], *p117*) in PA by two, formally equivalent, propositions[2].

---

[1] We assume that, in any system of Arithmetic such as IA, a function or relation containing free variables is recursive if and only if, for any given set of values for the free variables in the definition of the function or relation, the arithmetical value of the function, or the truth/falsity of the relation, can be determined in a finite number of steps from the *Axioms* of the system using its *Rules of Inference* by some mechanical procedure.

Another assumed characteristic of recursive expressions in a system of Arithmetic such as IA is that, for any given set of values for the free variables in the definition of the expression, it can be reduced, in a finite number of steps by some mechanical procedure, to an expression that consists of only a finite number of primitive symbols of a suitably constructed sub-language of IA such as PA, even though the definition and expression of the function or relation may involve an element of self-reference.

The first is a direct expression of the relation as the PA-formula $[q(\underline{k}, \underline{m})]$[3] - since a recursive arithmetical relation $q(x, y)$ can be expressed as a proposition in only the primitive symbols of PA once the variables $x$ and $y$ are replaced by the numerals $\underline{k}$ and $\underline{m}$[4], that represent the specific natural numbers $k$ and $m$.

The second is as a formally equivalent formula[5] $[Q(\underline{k}, \underline{m})]$ of PA that is defined in terms of the Gödel-Beta function by the usual Representation meta-Lemma ([Me64], *p131*). This meta-Lemma establishes that every recursive arithmetical relation, such as say $q(x, y)$, is instantiationally equivalent to an arithmetical relation $Q(x, y)$ that can be expressed in only the primitive symbols of PA.

Thus, although $q(x, y)$ may not necessarily be reducible to an expression that consists of only the primitive symbols of PA, $Q(x, y)$ is always a well-formed formula of PA. Further, $q(k, m)$ and $Q(k, m)$ are equivalent arithmetical propositions that can both be expressed in PA by the formally equivalent PA-propositions $[q(\underline{k}, \underline{m})]$ and $[Q(\underline{k}, \underline{m})]$, respectively, for any given natural numbers $k$ and $m$. It follows that though the relations $q(x, y)$ and $Q(x, y)$ are "arithmetically" equivalent, they are not "formally" equivalent.

I argue that this asymmetry yields alternative Representation and Self-reference meta-Lemmas that are critical to any exposition of Gödel's meta-reasoning. This reasoning can now be expressed avoiding any appeal to the truth of propositions of IA. I then argue that we can establish GUS as decidable in PA, and PA as omega-inconsistent. I argue further that Rosser's non-formal extension of Gödel's meta-reasoning is now seen as invalid.

---

[2] We use the terms "proposition" and "sentence" interchangeably. When referring to a well-formed symbolic expression, both terms imply that the expression has no free variables.

[3] We use square brackets to indicate that the expression (including square brackets) only denotes the PA-string that is named by the expression within the brackets. Thus, "$[q(\underline{k}, \underline{m})]$" is not part of the formal system PA. In this case, the PA-string named by "$q(k, m)$" is obtained by replacing the symbols constituting the arithmetical proposition "$q(k, m)$" by the PA-symbols of which they are the interpretations. The result is a PA-string of which "$q(k, m)$" is the arithmetic interpretation.

[4] We denote by "$\underline{n}$" the formal numeral in PA that represents the natural number "$n$" of IA.

[5] We use the formal term "formula" as corresponding to the intuitive term "expression". By "well-formed formula" we mean a symbolic expression that is constructed according to some grammatical rules of a system for the formation of symbol strings, by concatenation of the primitive symbols of the system, that are to be considered as "well-formed".

**1.2 The critical semantic elements of Gödel 's and Rosser's reasoning**

I argue that the critical elements in classical expositions (*such as* [Me64], *p145*) of Gödel's reasoning involve semantic interpretations in the Representation meta-Lemma, and in Gödel's Self-reference meta-Lemma. Both meta-Lemmas correlate the PA-provability of propositions with the truth of their interpretations in the standard model IA (*Intuitive Arithmetic*) of PA.

**1.3 The semantic form of the Representation meta-Lemma**

Classically (*eg.* [Me64], *p131, Proposition 3.23*), the Representation meta-Lemma essentially asserts, weakly and non-formally, the following:

> If $f(x_1, x_2, \ldots, x_n)$ is a recursive arithmetic relation of IA whose arguments are natural numbers, then there is another arithmetic relation $F(x_1, x_2, \ldots, x_n)$ of IA such that, for any given set of natural numbers $(k_1, k_2, \ldots, k_n, k_m)$ of IA:
>
> (*a*)  if $f(k_1, k_2, \ldots, k_n, k_m)$ is true in IA, then $[F(\underline{k}_1, \underline{k}_2, \ldots, \underline{k}_n, \underline{k}_m)]$ is PA-provable, and
>
> (*b*)  $[(E!x_m)F(\underline{k}_1, \underline{k}_2, \ldots, \underline{k}_n, x_m)]$ is PA-provable.

**1.4  The semantic form of Gödel's Self-reference meta-Lemma**

Classical expositions of Gödel's reasoning (*eg.* [Me64], *p145*) define a recursive arithmetic relation $q(x, y)$ that is true in IA if and only if $x$ is the Gödel-number of a well-formed formula [$H(z)$] of PA with a single free variable [$z$], and $y$ is the Gödel-number of a PA-proof of [$H(x)$].

So, in its classically semantic form, Gödel's Self-reference meta-Lemma essentially asserts, weakly and non-formally, the following:

> For any given formula [$H(z)$] of a single variable in PA whose Gödel-number is the natural number $h$ of IA, and any natural number $j$ of IA, we have, by the definition of $q(x, y)$, that:
>
> (c)    $q(h, j)$ is true in IA if, and only if, [$H(\underline{h})$] is PA-provable.

**1.5 A non-semantic Representation meta-Lemma**

In this paper, I replace the semantic Representation meta-Lemma §1.4 by the following stronger, meta-assertion:

> If $f(k_1, k_2, \ldots, k_n)$ is a recursive arithmetic relation of IA whose arguments are natural numbers then, by Hilbert and Bernays Representation meta-Lemma, there is another arithmetic relation $F(x_1, x_2, \ldots, x_n)$ of IA such that, for any given set of natural numbers $(k_1, k_2, \ldots, k_n)$ of IA:
>
> (d)    [$f(\underline{k}_1, \underline{k}_2, \ldots, \underline{k}_n) <=> F(\underline{k}_1, \underline{k}_2, \ldots, \underline{k}_n)$] is PA-provable.

**1.6 A non-semantic Self-reference meta-Lemma**

Similarly, I replace the semantic Self-reference meta-Lemma §1.5 by the following stronger, meta-assertions:

> For any given formula [$H(z)$] of a single variable in PA whose Gödel-number is the natural number $h$ of IA, and any natural number $j$ of IA, we have, by the definition of $q(x, y)$, that:

(*e*)    [$q(\underline{h}, \underline{j}) \Rightarrow H(\underline{h})$] is PA-provable.

(*f*)    If [$H(\underline{h})$] is PA-provable, then [$q(\underline{h}, \underline{k})$] is PA-provable for some natural number *k* of IA.

### 1.7 Consequences

I argue that §1.6(*d*) and §1.7(*e*) are denumerable sets of non-semantic meta-Lemmas in PA, whilst §1.7(*f*) is a non-semantic meta-Lemma in the meta-theory of PA.

I also argue that Goedel's reasoning can be constructively interpreted as implying that, from the PA-provability of [$(Ex)H(x)$], we may not always assume the existence of some numeral [$\underline{h}$] such that [$H(\underline{h})$] is provable in PA. All we may conclude is that ~$H(x)$ may not be a uniformly decidable predicate in IA; however, ~$H(n)$ may yet hold in IA for any given natural number *n* - which may be provable meta-mathematically. I argue that such an invalid assumption is, however, implicit in Rosser's proof of undecidability.

I argue that, by using the meta-Lemmas §1.6 and §1.7 to recast classical expositions of Gödel's and Rosser's semantic proofs of undecidability, we can determine that PA is omega-inconsistent and that Rosser's assumption (*as above*) is invalid.

Consequently, Gödel's and Rosser's semantic arguments do not yield formally undecidable propositions in PA in a constructive and intuitionistically unobjectionable way.

### 1.8 Background and conclusions

The above is an extension of earlier arguments questioning various disquieting features of the classical expositions and interpretation of Gödel's reasoning, such as in Mendelson [Me64].

In Anand [An01], I question the significance, and relevance, of standard PA. I argue the thesis that differing perceptions of the interpretation of Gödel's meta-reasoning, and conclusions, arise because the system of standard first order Peano Arithmetic PA, in which Gödel's meta-reasoning is classically considered, is only one - and perhaps not the most representative - of several, significantly differing, systems ([An01], §3.0) that can be defined to formalise our system of *Intuitive Arithmetic* IA of the natural numbers (*which we take both as the Arithmetic based on an intuitive interpretation of Dedekind's formulation of Peano's Postulates, and as the standard interpretation of any formalisation of these Postulates*).

In Anand [An02], I outline in general terms a constructive, and intuitionistically unobjectionable, formalisation PP of our *Intuitive Arithmetic* IA of the natural numbers in which the *Axioms* and *Rules of Inference* are recursively definable. I argue that, although Gödel's undecidable sentence GUS is a well-formed formula in PP ([An02], *§3.9*), it is an ill-defined proposition that formally reflects the Liar paradox in PP ([An02], *§4.5*).

In this paper I assume some familiarity with the issues addressed in the above two papers, especially those concerning the asymmetry ([An02], *§6.1*) between a primitive recursive function of IA, and the interpretation in IA of its strong representation ([An01], *§2.8*) in PA. I argue that this asymmetry can also be expressed formally, and use this to develop a formal expression of Gödel's meta-reasoning to argue that:

(*i*)   Gödel's undecidable sentence GUS is actually decidable in PA under a reasonable interpretation of his meta-reasoning that, avoiding any appeal to the truth of propositions of IA, is constructive and intuitionistically unobjectionable;

(*ii*)  PA is not omega-consistent under such interpretation;

(*iii*) If we take IA to be the standard model for PA, this implies that standard PA is semantically inconsistent under the standard interpretations of Gödel's meta-reasoning;

(*iv*) Rosser's extension of Gödel's meta-reasoning, establishing an undecidable proposition RUS in a consistent PA, is invalid under such an interpretation.

## 2. The formal system PA

### 2.1. Peano's Arithmetic PA

We start by noting that, in standard (*first order*) *Peano's Arithmetic* PA, the following *Axioms* and *Rules of Inference* (*essentially* [Me64], *p103*) can reasonably be taken as our (*sole*) basis for assigning provability values to the well-formed sentences of PA:

(*PA1*)  $\sim(\underline{0} = (x+\underline{1}))$

(*PA2*)  $\sim(x = y) \Rightarrow \sim((x+\underline{1}) = (y+\underline{1}))$

(*PA3*)  $(x+\underline{0}) = x$

(*PA4*)  $(x+(y+\underline{1})) = ((x+y)+\underline{1})$

(*PA5*)  $(x*\underline{0}) = \underline{0}$

(*PA6*)  $(x*(y+\underline{1})) = ((x*y)+x)$

### 2.2. The alphabet S of PA

These formulas and assertions are expressed using only a small set S of undefined, primitive, symbols such as "+" (*addition*) and "*" (*multiplication*) only apart from the logical symbols "~"(*not*), "=>" (*implies*), "=" (*equals*), "&" (*and*), "v" (*or*), "[(A*x*)]" (*a special symbol for representing 'all x'*), "[(E*x*)]" (*there exists x*), "*x, y, ...*" (*variables*),

"*a, b, c, ...*" (*constants*), "<u>0</u>" (*numeral zero*), "<u>1</u>" (*numeral one*), and the two parentheses "(", ")".

The non-terminating series of numerals "<u>0</u>, <u>0</u>+<u>1</u>, (<u>0</u>+<u>1</u>)+<u>1</u>, ((<u>0</u>+<u>1</u>)+<u>1</u>)+<u>1</u>, ...", is taken as formally representing, in PA, the various non-terminating, intuitive, natural number series such as "0, 1, 10, 11, ..." (*in binary format*), or "0, 1, 2, 3, ..." (*in the more common decimal format*).

### 2.3. *Rules of Inference* in PA

We take as our (*only*) means for deriving other provable assertions in PA, from the basic *Axioms* (*PA1-6*), the standard first order logical axioms ([Me64], *p57*), and the following *Rules of Inference*, namely *Modus Pon*ens, *Induction* and *Generalisation*:

- (*PAR1*)   *Modus Ponens*: From [*F*] and [*F* => *G*] we may conclude [*G*], where [*F*] and [*G*] are any well-defined formulas of PA.

- (*PAR2*)   *Induction*: From [*F*(<u>0</u>)] and [(A*x*)(*F*(*x*) => *F*(*x*+<u>1</u>))] we may conclude [(A*x*)*F*(*x*)].

- (*PAR3*)   *Generalisation*: From [*F*(*x*)] we may conclude [(A*x*)*F*(*x*)].

## 3. The significance of the formal system PA

### 3.1. Meta-assertions about PA can be expressed in PA

We note that, as described by Gödel, terms such as "well-formed formula", "well-formed proposition" and "proof sequence" can be formally defined in the language of PA. We can thus express meta-assertions about PA such as " '[*F*(*x*)]' is a well-formed formula in PA", " '[(A*x*)*F*(*x*)]' is a well-formed proposition in PA", " '[*F*(*x*)]' is a provable formula in PA", " '[(A*x*)*F*(*x*)]' is a provable proposition in PA" and "The 'Gödel' proposition is not provable in PA", amongst others, as equivalent arithmetical assertions in PA (*Gödel 1931, pp.14*).

## 3.2. Provability in PA

We note that Gödel defines a well-formed formula *P* of IA as provable in PA if and only if there is a finite proof sequence, consisting of well-formed formulas of PA, each of which is either an *Axiom* (*PA1-6*), or an immediate consequence of the preceding provable well-formed formulas by the *Rules of Inference* (*PAR1-3*), and where *P* is the final well-formed formula in the sequence. The *Axioms* (*PA1-6*) are, of course, all provable well-formed formulas in PA by definition ([Go31], *p13*).

The provable propositions in PA are characterised by the fact that they are all provable well-formed formulas without variables, such as "$[(Ax)F(x)]$", expressed using only the small set of primitive, undefined symbols of the alphabet S.

## 3.3. Recursive functions of IA are not directly expressible in PA

The significance of the alphabet S selected for expressing provable propositions of IA is that many significant arithmetical functions of IA such as "*n!*" (*factorial*), "*m^n*" (*exponential*), "*m/n*" (*division*), "*n* is a prime number", amongst others, are defined recursively, and recursive functions are Turing-computable ([Me64], *p237*).

Thus we note that *n!*, for instance, is defined by:

(*i*)   $0! = 1$

(*ii*)   $(n+1)! = (n+1)*(n!)$ for all natural numbers *n*.

Now, for a given natural number *n*, any expression involving *n!* can clearly be reduced in a finite number of steps to an expression that consists of only the symbols of the alphabet

S, eg. "3! = 5" to "3*(2*(1*(1))) = 5". Thus, in this case, "3! = 5" can be viewed simply as a shorthand notation (*i.e. another name, or definition*) for the assertion "3*(2*(1*(1))) = 5".

However, if we attempt to eliminate the symbol "!" on the right hand side of (*ii*) for an unspecified $x$, we soon discover that "$x$!" is not directly reducible to a form that is expressible in only the symbols of the alphabet S.

The question thus arises: Can we locate a general expression FACT($x$) of IA that is expressible in only the alphabet S, and is such that its values are identical to those of $x$! whenever a natural number $n$ is substituted for $x$?

**3.4. A non-semantic Representation meta-Lemma**

(*a*) This issue is addressed, and answered indirectly[6], by a meta-Lemma of Hilbert and Bernays, which establishes that every recursive arithmetical relation definable in

---

[6] Indirectly in the sense that, for a recursive arithmetical function such as say "$x$!", where $x$ ranges over the natural numbers we can, loosely speaking, define a denumerable series of recursive arithmetical functions FACT$_n$($x$) such that ([Me64], *p118 & 131*):

(*i*) FACT$_n$($x$) can be strongly represented in PA for any natural number $n$;

(*ii*) For every natural number $n$:

FACT$_n$($x$) = $x$! for all natural numbers $0 =< x =< n$.

In other words, for a given natural number $n$, FACT$_n$($x$) is a recursive arithmetical function, termed the Gödel-Beta function that, loosely speaking, generates the finite sequence 0!, 1!, 2!, ... , $n$!.

We note, however, that if $n<m$, then the value of FACT$_n$($m$) is not equal to $m$!.

We also note that Gödel's reasoning in the proof of his undecidability theorem essentially concerns itself not with an original recursive arithmetical function such as $x$! but, loosely speaking, with the provability in PA of the strong representation of the associated Gödel-Beta functions FACT$_n$($x$). I attempt to highlight this factor in Anand, 2001, §2_9.

I argue that it is this asymmetry, at the core of Gödel's reasoning, that is actually reflected in the curious

IA, say $f(x) = y$ where $f(x)$ is a recursive arithmetical function of IA, is indeed instantiationally equivalent to another, uniquely defined, arithmetical relation [$F(x, y)$] of IA that can be represented in PA (*i.e. F(x, y) is a well-formed formula of IA that is expressible in only the alphabet S*), such that (*cf.* [Me64], *p118-120*):

(*i*)     $(Ak)(Am)$PA proves: $[(f(\underline{k}) = \underline{m}) <=> F(\underline{k}, \underline{m})]$

(*ii*)     PA proves: $[(E!y)F(\underline{k}, y)]$ [7]

where [$\underline{n}$] denotes the numeral in PA that represents the natural number $n$ of IA, and "E!" denotes uniqueness of the existential assertion.

(*b*) We note that (*a*)(*i*) expresses a denumerable set of Lemmas in PA since, for given natural numbers $k$ and $m$, any recursive arithmetical proposition such as "$f(k) = m$" can be expressed as a proposition "$f(\underline{k}) = \underline{m}$" in terms of the alphabet S only. It is thus a valid well-formed formula of PA that is decidable[8] in PA, since $f(k)$ is a Turing-computable function of IA ([Me64], *p231*).

For instance, if $f(n)$ is the factorial function $n!$, then, as mentioned above in §3.3, the proposition "$3! = 5$" of IA can also be expressed as "$3*(2*(1*(1))) = 5$", which is the interpretation in IA of the decidable PA-proposition $[\underline{3}*(\underline{2}*(\underline{1}*(\underline{1}))) = \underline{5}]$.

---

"true but unprovable" elements in Gödel's conclusions. Gödel, however, noted that his conclusions were a reflection of implicit, transfinite, elements of PA ([Go31], *footnote 48a*).

[7] We note the asymmetry that, since $f(x)$ may not, as in the case of $x!$, be expressible in the alphabet S of PA, we do not correspondingly have that PA proves: $[(E!y)(f(\underline{k}) = y)]$ for every recursive arithmetical function $f(x)$ of IA, since $(E!y)(f(\underline{k}) = y)$ may not be the interpretation of a well-formed formula of PA.

[8] Gödel 1931, pp.23 footnote 39

**3.5. Gödel-numbering**

We note next that every expression [F], constructed by concatenation from the primitive, undefined symbols of S, can be assigned a unique natural number of IA, which we term as the "Gödel-number" of the expression [F] ([Go31], *p13*).

Now Gödel has established ([Go31], *p17-22*) that we can define a recursive arithmetic relation *prf(x, y)* in IA (*Gödel's 'yBx'*), constructed out of 44 "simpler" recursive arithmetic functions and relations, such that, for any given natural numbers $k$ and $m$, *prf($\underline{k}$, $\underline{m}$)* is true in IA if and only if $m$ is the Gödel-number of a finite proof sequence *PRF* in PA for some well-formed formula [K] in PA whose Gödel-number is $k$.

However, as we argue in §3.4(b), for any given natural numbers $k$ and $m$, [*prf($\underline{k}$, $\underline{m}$)*] is also provable in PA if and only if $m$ is the Gödel-number of a finite proof sequence *PRF* in PA for some well-formed formula [K] in PA whose Gödel-number is $k$.

**3.6. Gödel's formal Self-reference Lemmas**

From the recursive arithmetic relation *prf(x, y)*, we can thus define another recursive arithmetic relation *q(x, y)* (*Mendelson's $W_1(x, y)$,*[Me64], *p143*) in IA which is similarly true in IA if and only if $x$ is the Gödel-number of a well-formed formula [H(z)] of PA with a single free variable [z], and $y$ is the Gödel-number of a proof of [H(x)] in PA.

Hence the constructive self-reference, which lies at the core of Gödel's meta-reasoning, is that, since [*q(h, j)*] is a PA-formula for any given natural numbers $h, j$, it is provable in PA if and only if the expression *J* in PA, whose Gödel number is $j$, is a proof sequence in PA of the well-formed proposition [H($\underline{h}$)] of PA, where $h$ is the Gödel-number of the formula [H(z)] in PA containing only one free variable. More precisely, if [H(z)] is as defined above, we have that, by definition:

(i)   (A*j*)PA proves: [*q*(*h̲*, *j̲*) => *H*(*h̲*)]

(ii)   PA proves: [*H*(*h̲*)] => (E*j*)PA proves: [*q*(*h̲*, *j̲*)]

We note, firstly, that whilst (*i*) is a denumerable set of Lemmas in PA, (*ii*) is actually a meta-Lemma in the meta-theory of PA. Further that, here also, as argued in §3.4(*b*), (*i*) and (*ii*) hold since, for given natural numbers *h* and *j*, any recursive arithmetical expression such as *q*(*h*, *j*) can also be expressed directly in terms of the alphabet S only. It can thus be shown to be a valid, decidable, well-formed formula of PA.

## 4. Gödel's meta-reasoning

### 4.1. The two basic Lemmas and GUS

Now, by the Representation meta-Lemmas, the recursive arithmetic relation *q*(*x*, *y*) is also instantiationally equivalent to another, uniquely defined, well-formed arithmetic relation *Q*(*x*, *y*), such that:

(i)   (A*k*)(A*m*)PA proves: [*q*(*k̲*, *m̲*) <=> *Q*(*k̲*, *m̲*)]

If *p* is the Gödel-number of the well-formed formula [(A*y*)(~*Q*(*x*, *y*))], we consider the well-formed "Gödelian" proposition GUS expressed by [(A*y*)(~*Q*(*p̲*, *y*))].

By Gödel's Self-reference meta-Lemmas, we then have that:

(ii)   (A*j*)PA proves: [*q*(*p̲*, *j̲*) => (A*y*)(~*Q*(*p̲*, *y*))]

(iii)   PA proves: [(A*y*)(~*Q*(*p̲*, *y*))] => (E*j*)PA proves: [*q*(*p̲*, *j̲*)]

### 4.2. Gödel's semantic meta-mathematical proof of undecidability[9]

(*a*) We assume firstly that *r* is the Gödel-number of some proof-sequence *R* in PA for the proposition [(A*y*)(~*Q*(*p*, *y*))]. It then follows from Gödel's Self-reference lemma §4.1(*iii*) that *q*(*p*, *r*) is true in IA. By the Representation Lemmas §4.1(*i*), this implies that [*Q*(*p*, *r*)] is provable in PA. However, assuming standard logical axioms for PA, from the provability of [(A*y*)(~*Q*(*p*, *y*))] in PA, we have that [~*Q*(*p*, *r*)] is provable in PA. It follows that there is no natural number *r* that is the Gödel-number of a proof-sequence *R* in PA for the proposition [(A*y*)(~*Q*(*p*, *y*))]. Hence, by §4.1(*ii*), *q*(*p*, *r*) is not true in IA for any *r*, and so [(A*y*)(~*Q*(*p*, *y*))] is not provable in PA.

(*b*) We assume next that *r* is the Gödel-number of some proof-sequence *R* in PA for the proposition [~(A*y*)(~*Q*(*p*, *y*))]. However, we have by §4.2(*a*) that *q*(*p*, *r*) is not true in IA for any *r*, and so, by §4.1(*i*), [~*Q*(*p*, *r*)] is provable in PA for all *r*. Assuming that PA is omega-consistent, it follows that [~(A*y*)(~*Q*(*p*, *y*))] is not provable in PA.

( *Note: We define PA as omega-consistent if and only if there is no formula '[F(x)]' such that '(n)PA proves: [F(n)]' and 'PA proves: [~(Ax)F(x)]'.*)

(*c*) The classical conclusion drawn from the above is that, if PA is omega-consistent, then the proposition [(A*y*)(~*Q*(*p*, *y*))] is undecidable in PA, since neither [(A*y*)(~*Q*(*p*, *y*))] nor [~(A*y*)(~*Q*(*p*, *y*))] is provable in PA ([Me64], *pp.143-144*).

(*d*) By §4.2(*a*), we have that [~*Q*(*p*, *r*)] is provable in PA for all natural numbers *r*. We thus have that ~*Q*(*p*, *y*) is true in IA for all *y* (*since the domain of IA is expressible*

---

[9] See also Anand, 2001, §5.8 and Anand, 2002, §3.9.

*in PA, and every provable proposition of PA is true in IA*). It now follows by our definition of quantification in IA (*ignoring this point lends validity to Lucas' Gödelian argument*) that "$(Ay)(\sim Q(\underline{p}, y))$" can be asserted as a true proposition in IA.

(*e*) It is only in this sense that, although $[(Ay)(\sim Q(\underline{p}, y))]$ is not provable in PA, Gödel was able to assert "$(Ay)(\sim Q(\underline{p}, y))$" as a true proposition in IA. However, I argue here that the truth of the latter proposition in IA is clearly of a definitional nature, and a formal consequence of the *Axioms* and *Rules of Inference* of PA.

(*f*) I argue that the curious conclusion (*e*) reflects the fact that, in the presence of a *Generalisation Rule of Inference*, we cannot correspondingly postulate (*or define*) formally that $[(Ay)(\sim Q(\underline{p}, y))]$ follows as a provable proposition in PA from the provability of $[\sim Q(\underline{p}, \underline{r})]$ in PA for all natural numbers $r$, if PA is assumed consistent [*Anand 2001, §1.11*].

**4.3. A non-semantic expression of Gödel's meta-proof of undecidability**

(*a*) Now we note that we can also express the above meta-reasoning non-semantically as the deduction:

(*i*)      PA proves: $[(Ay)(\sim Q(\underline{p}, y))]$

        ... *Hypothesis*

(*ii*)     PA proves: $[(Ay)(\sim Q(\underline{p}, y))] \Rightarrow (Er)$PA proves: $[q(\underline{p}, \underline{r})]$

        ... *By the Self-reference meta-lemmas §4.1(iii)*

(*iii*)    $(Er)$PA proves: $[q(\underline{p}, \underline{r})]$

        ... *From (i) and (ii) by Modus Ponens in the meta-theory*

 (*iv*)  (A*r*)PA proves: [*q*(*p*, *r*) <=> *Q*(*p*, *r*)]

    *... By the Representation lemmas §4.1(i)*

 (*v*)  (E*r*)PA proves: [*Q*(*p*, *r*)]

    *... From (iii) and (iv) by Modus Ponens, assuming (i)*

 (*vi*)  PA proves: [(E*y*)*Q*(*p*, *y*)]

    *... From (v) by the logical axioms of PA, assuming (i)*

 (*vii*)  PA proves: [(A*y*)(~*Q*(*p*, *y*)) => (E*y*)(~*Q*(*p*, *y*))]

    *... Tautology*

 (*viii*)  PA proves: [(E*y*)(~*Q*(*p*, *y*))]

    *... From (i) and (vii) by Modus Ponens*

Curiously, avoiding the inference that "PA proves: [(*i*) => (*vi*)]" (*by the Deduction Theorem*), the classical conclusion drawn from the above is that, since (*vi*) contradicts (*viii*), [(A*y*)(~*Q*(*p*, *y*))] is not provable in a consistent PA, i.e.:

 (*ix*)  ~PA proves: [(A*y*)(~*Q*(*p*, *y*))]

    *... assuming PA is consistent*

(*b*) We also have the further classical deduction:

 (*i*)  PA proves: [~(A*y*)(~*Q*(*p*, *y*))]

    *... Hypothesis*

 (*ii*)  PA proves: [(E*y*)*Q*(*p*, *y*)]

       *... From (i) by definition*

(*iii*)    (A*r*)PA proves: [*q*(*p*, *r*) => (A*y*)(~*Q*(*p*, *y*)]

       *... By the Self-reference lemmas §4.1(ii)*

(*iv*)    (A*r*)PA proves: [~(A*y*)(~*Q*(*p*, *y*) => ~*q*(*p*, *r*)]

       *... From (iii) by the logical axioms of PA*

(*v*)    (A*r*)PA proves: [~*q*(*p*, *r*)]

       *... From (i) and (iv) by Modus Ponens*

(*vi*)    (A*r*)PA proves: [*q*(*p*, *r*) <=> *Q*(*p*, *r*)]

       *... By the Representation lemmas §4.1(i)*

(*vii*)    (A*r*)PA proves: [~*q*(*p*, *r*) <=> ~*Q*(*p*, *r*)]

       *... From (vi) by the logical axioms of PA*

(*viii*)    (A*r*)PA proves: [~*Q*(*p*, *r*)]

       *... From (v) and (vii) by Modus Ponens, assuming (i)*

Since (*i*) and (*viii*) contradict the omega-consistency of PA, the further conclusion drawn classically from the above reasoning is that [~(A*y*)(~*Q*(*p*, *y*))] too is not provable in an omega-consistent PA, i.e.:

(*ix*)    ~PA proves: [~(A*y*)(~*Q*(*p*, *y*))]

   ... *Assuming PA is omega-consistent.*

So we have the classical conclusion drawn from Gödel's above meta-reasoning that GUS is an undecidable proposition in an omega-consistent PA.

## 5. Can GUS be decidable, and PA omega-inconsistent?

However, as noted in §4.3(*a*)(*viii*), it can also be argued reasonably in the meta-theory of PA, assuming *p* is the Gödel-number of the well-formed formula [(A*y*)(~*Q*(*x*, *y*))], that:

(*i*)   (A*r*)PA proves: [*q*(*p*, *r*) => (A*y*)(~*Q*(*p*, *y*))]

   ... *By the Self-reference lemmas §4.1(ii)*

(*ii*)  PA proves: [(A*y*)(~*Q*(*p*, *y*))] => (E*r*)PA proves: [*q*(*p*, *r*)]

   ... *By the Self-reference meta-lemmas §4.1(iii)*

(*iii*) (A*r*)PA proves: [*q*(*p*, *r*) <=> *Q*(*p*, *r*)]

   ... *By the Representation lemmas §4.1(i)*

(*iv*)  PA proves: [(A*y*)(~*Q*(*p*, *y*))] => (E*r*)PA proves: [*Q*(*p*, *r*)]

   ... *From (ii) and (iii)*

(*v*)   PA proves: [(A*y*)(~*Q*(*p*, *y*))] => PA proves: [(E*y*)*Q*(*p*, *y*)]

   ... *From (iv) by the logical axioms of PA*

(*vi*)  PA proves: [(A*y*)(~*Q*(*p*, *y*)) => (E*y*)*Q*(*p*, *y*)]

   *... From (v) by the Deduction theorem*

(*vii*)  PA proves: [(A*y*)(~*Q*(*p*, *y*)) => ~(A*y*)(~*Q*(*p*, *y*))]

   *... From (vi) by definition*

(*viii*) PA proves: [~(A*y*)(~*Q*(*p*, *y*))]

   *... From (vii) assuming PA is consistent*

We can thus alternatively, and prima facie quite reasonably, conclude from Gödel's meta-reasoning that, if PA is consistent, then ~(A*y*)(~*Q*(*p*, *y*)) is a theorem of PA, and not an undecidable proposition in PA. Now we also have that:

(*ix*)  PA proves: [(E*y*)*Q*(*p*, *y*)]

   *... From (viii) by definition*

(*x*)  (A*r*)PA proves: [~(A*y*)(~*Q*(*p*, *y*) => ~*q*(*p*, *r*)]

   *... From (i) by the logical axioms of PA*

(*xi*)  (A*r*)PA proves: [~*q*(*p*, *r*)]

   *... From (viii) and (x) by Modus Ponens*

(*xii*) (A*r*)PA proves: [~*q*(*p*, *r*) <=> ~*Q*(*p*, *r*)]

   *... By the Representation lemmas §4.1(i)*

  (*xiii*) (A*r*)PA proves: [~$Q(\underline{p}, \underline{r})$]

    *... From (xi) and (xii) by Modus Ponens*

Thus, from (*viii*) and (*xiii*), we may further conclude that a consistent PA is not omega-consistent.

## 6. Standard PA is semantically inconsistent

On this interpretation of Gödel's meta-reasoning, it follows that, in every model M of PA, (E*y*)$Q(p^*, y)$ is a true proposition in M, even though $Q(p^*, r^*)$ would be false for every $r^*$ in M that is the interpretation of a numeral $\underline{r}$ of PA.

If IA is taken to be the classical standard interpretation of PA, then both ~(A*y*)(~$Q(p, y)$) and (A*y*)(~$Q(p, y)$) interpret as true propositions in IA, and the above meta-reasoning then implies that standard PA is semantically inconsistent[10].

## 7. Is Rosser's undecidable proposition established constructively?

The question now arises: Are other undecidable propositions also formally derivable in PA on the above interpretation of Gödel's meta-reasoning?

We address only one instance of this question by considering Mendelson's version ([Me64], *p144-146*) of Rosser's extension of Gödel's meta-reasoning, where he seeks to establish an undecidable proposition RUS in a consistent, standard PA.

---

[10] It might be argued that if (E*y*)$Q(p, y)$ is a true assertion in IA, then $Q(p, k)$ must be true in IA for some natural number $k$. Since $Q(x, y)$ is recursive, it would then follow that [$Q(p^*, k)$] is PA-provable. Hence the above semantical inconsistency of IA implies the formal inconsistency of PA.

However we can reasonably argue against such reasoning and, as we do in the following review of Rosser's proof, hold that we may not conclude the constructive existence of the natural number $k$ from the PA-provability of [(E*y*)$Q(p^*, y)$].

## 7.1. Rosser's Self-reference Lemmas

Rosser reasoning is that, along with the recursive arithmetic relation $q(x, y)$, we can also define a recursive arithmetical relation $s(x, y)$ in IA which is provable in PA if and only if $x$ is the Gödel-number of a well-formed formula $[H(z)]$ of PA with a single free variable $[z]$, and $y$ is the Gödel-number of a proof of $[\sim H(x)]$ in PA.

Since $[s(h, j)]$ too is expressible in PA for given natural numbers $h, j$, it also is provable in PA if and only if the expression $J$ in PA, whose Gödel number is $j$, is a proof sequence in PA of the well-formed proposition $[\sim H(\underline{h})]$ of PA, where $h$ is the Gödel-number of the formula $[H(z)]$ in PA.

More precisely, if $[H(z)]$ is as defined above, we have that, by definition:

(*i*)   (A*j*)PA proves: $[s(\underline{h}, \underline{j}) \Rightarrow \sim H(\underline{h})]$

(*ii*)   PA proves: $[\sim H(\underline{h})] \Rightarrow$ (E*j*)PA proves: $[s(\underline{h}, \underline{j})]$

## 7.2. Rosser's undecidable proposition RUS

Now, by Hilbert and Bernays Representation meta-Lemma, the recursive relation $s(x, y)$ too is instantiationally equivalent to another, uniquely defined, well-formed arithmetical relation $S(x, y)$, so that:

(*i*)   (A*k*)(A*m*)PA proves: $[s(\underline{k}, \underline{m}) \Leftrightarrow S(\underline{k}, \underline{m})]$

If $u$ is the Gödel-number of the well-formed formula $[(Ay)(Q(x, y) \Rightarrow (Ez)(z =< y \,\&\, S(x, z)))]$ [11], we consider then the well-formed Rosser proposition RUS expressed by $[(Ay)(Q(\underline{u}, y) \Rightarrow (Ez)(z =< y \,\&\, S(\underline{u}, z)))]$.

---

[11] We use the symbolism "=<" to denote "equal to or less than".

Clearly, by definition, we have the denumerable set of Rosser's non-semantic Self-reference meta-Lemmas:

(*ii*)     (A*j*)PA proves: $[q(\underline{u}, \underline{j}) \Rightarrow (Ay)(Q(\underline{u}, y) \Rightarrow (Ez)(z=<y \& S(\underline{u}, z)))]$

(*iii*)     PA proves: $[(Ay)(Q(\underline{u}, y) \Rightarrow (Ez)(z=<y \& S(\underline{u}, z)))] \Rightarrow (Ej)$PA proves: $[q(\underline{u}, \underline{j})]$

(*iv*)     (A*j*)PA proves: $[s(\underline{u}, \underline{j}) \Rightarrow \sim(Ay)(Q(\underline{u}, y) \Rightarrow (Ez)(z=<y \& S(\underline{u}, z)))]$

(*v*)     PA proves: $[\sim(Ay)(Q(\underline{u}, y) \Rightarrow (Ez)(z=<y \& S(\underline{u}, z)))] \Rightarrow (Ej)$PA proves: $[s(\underline{u}, \underline{j})]$

**7.3. Rosser's semantic meta-mathematical proof of undecidability**

(*a*) We assume now that *r* is the Gödel-number of some proof sequence *R* in PA for the proposition $[(Ay)(Q(\underline{u}, y) \Rightarrow (Ez)(z=<y \& S(\underline{u}, z)))]$. It then follows from Rosser's Self-reference meta-Lemma §7.2(*iii*) that $q(u, r)$ is true in IA. Also, by the Representation meta-Lemmas §4.1(*i*), this implies that $[Q(\underline{u}, \underline{r})]$ is provable in PA. However, assuming standard logical axioms for PA, from the provability of $[(Ay)(Q(\underline{u}, y) \Rightarrow (Ez)(z=<y \& S(\underline{u}, z)))]$ in PA, we have that $[Q(\underline{u}, \underline{r}) \Rightarrow (Ez)(z=<\underline{r} \& S(\underline{u}, z))]$ is provable in PA and, by Modus Ponens, that $[(Ez)(z=<\underline{r} \& S(\underline{u}, z))]$ is provable in PA. Now, since PA is consistent, and $[(Ay)(Q(\underline{u}, y) \Rightarrow (Ez)(z=<y \& S(\underline{u}, z)))]$ is provable in PA, it follows there is no proof in PA of $[\sim(Ay)(Q(\underline{u}, y) \Rightarrow (Ez)(z=<y \& S(\underline{u}, z)))]$. So, by §7.2(*iv*), $s(u, n)$ is not true in IA for any natural number *n*, and so $\sim s(u, n)$ is true in IA for every natural number *n*. Hence, by §7.2(*i*), $[\sim S(\underline{u}, \underline{n})]$ too is provable in PA for every natural number *n*. We thus have, by the standard logical axioms of PA, that $[\sim(Ez)(z=<\underline{r} \& S(\underline{u}, z))]$ is provable in PA. Since this contradicts our earlier derivation, it follows that $[(Ay)(Q(\underline{u}, y) \Rightarrow (Ez)(z=<y \& S(\underline{u}, z)))]$ is not provable in a consistent PA.

(*b*) We assume next that *r* is the Gödel-number of some proof-sequence *R* in PA for the proposition [~(A*y*)(*Q*(*u*, *y*) => (E*z*)(*z*=<*y* & *S*(*u*, *z*)))]. By §7.2(*v*), we then have that *s*(*u*, *r*) is true in IA. Hence, by §7.2(*i*), [*S*(*u*, *r*)] is provable in PA. However, if PA is consistent, we have by §7.3(*a*) that there is no proof sequence in PA for [(A*y*)(*Q*(*u*, *y*) => (E*z*)(*z*=<*y* & *S*(*u*, *z*)))], and so, by §7.2(*iii*), ~*q*(*u*, *n*) is true in IA for every natural number *n*. Hence, by §4.1(*i*), [~*Q*(*u*, *n*)] is provable in PA for all natural numbers *n*.

   (*i*) We thus have, by the standard logical axioms of PA, that [*y*=<*r* => ~*Q*(*u*, *y*)] is provable in PA.

   We consider now the following deduction:

      (*1*)   [*r*=<*y*]

           ... *Hypothesis*

      (*2*)   [*S*(*u*, *r*)]

           ... *By (b) above*

      (*3*)   [*r*=<*y* & *S*(*u*, *r*)]

           ... *From (1) and (2), Tautology*

      (*4*)   [(E*z*)(*z*=<*y* & *S*(*u*, *z*))]

           ... *From (3) by the logical axioms of PA*

   (*ii*) From (*1*)-(*4*), by the Deduction Theorem, we have that [*r*=<*y* => (E*z*)(*z*=<*y* & *S*(*u*, *z*))] is provable in PA.

    *(iii)* Now, by the logical axioms of PA, we have that $[y=\underline{r} \vee \underline{r}=\underline{y}]$ is provable in PA.

    *(iv)* From *(i)*-*(iii)*, also by the logical axioms of PA, we now have that $[\sim Q(\underline{u}, y) \vee (Ez)(z=\underline{y} \& S(\underline{u}, z))]$ is provable in PA.

    *(v)* We have further, again by the logical axioms of PA, that $[(Ay)(Q(\underline{u}, y) => (Ez)(z=\underline{y} \& S(\underline{u}, z)))]$ is provable in PA.

However, this contradicts our assumption that the proposition $[\sim(Ay)(Q(\underline{u}, y) => (Ez)(z=\underline{y} \& S(\underline{u}, z)))]$ is provable in PA.

The classical conclusion from the above meta-reasoning is that $[\sim(Ay)(Q(\underline{u}, y) => (Ez)(z=\underline{y} \& S(\underline{u}, z)))]$ too is not provable in a consistent PA, and so RUS is undecidable in PA ([Me64], *p146*).

**7.4. A non-semantic expression of Rosser's meta-proof of undecidability**

However, we note that the Rosser's above, semantic, meta-reasoning contains an invalid deduction that is exposed when the argument is expressed non-semantically as below.

*(a)* Now, the first half of the argument can be expressed as:

    *(i)* PA proves: $[(Ay)(Q(\underline{u}, y) => (Ez)(z=\underline{y} \& S(\underline{u}, z)))]$

        ... *Hypothesis*

    *(ii)* PA proves: $[(Ay)(Q(\underline{u}, y) => (Ez)(z=\underline{y} \& S(\underline{u}, z)))] => (Er)$PA proves: $[q(\underline{u}, \underline{r})]$

        ... *By the Self-reference meta-Lemma §7.2(iii)*

(*iii*)    (E*r*)PA proves: [*q*(*u*, *r*)]

    *... From (i) and (ii) by Modus Ponens in the meta-theory*

(*iv*)    (A*r*)PA proves: [*q*(*u*, *r*) <=> *Q*(*u*, *r*)]

    *... By the Representation meta-Lemmas §4.1(i)*

(*v*)    (E*r*)PA proves: [*Q*(*u*, *r*)]

    *... From (iii) and (iv) by Modus Ponens, assuming (i)*

(*vi*)    (A*r*)PA proves: [*Q*(*u*, *r*) => (E*z*)(*z*=<*r* & *S*(*u*, *z*))]

    *... From (i) by the logical axioms of PA*

(*vii*)    (E*r*)PA proves: [(E*z*)(*z*=<*r* & *S*(*u*, *z*))]

    *... From (v) and (vi) by Modus Ponens, assuming (i)*

(*viii*)    ~PA proves: [~(A*y*)(*Q*(*u*, *y*) => (E*z*)(*z*=<*y* & *S*(*u*, *z*)))]

    *... From (i), assuming PA is consistent*

(*ix*)    (A*r*)PA proves: [~*s*(*u*, *r*))]

    *... From (i) by the Self-reference meta-Lemmas §7.2(iv)*

(*x*)    (A*r*)PA proves: [~*S*(*u*, *r*)]

    *... From (ix) by the Representation meta-Lemmas §7.2(i), assuming (i)*

(*xi*)    (A*r*)PA proves: [(A*z*)(*z*=<*r* => ~*S*(*u*, *z*))]

    *... From (x) by the logical axioms of PA, assuming (i)*

(*xii*)   (A*r*)PA proves: [~(E*z*)(*z*=<*r* & *S*(*u*, *z*))]

   ... *From (xi) by the logical axioms of PA, assuming (i)*

Since (*xii*) contradicts (*vii*), the classical conclusion from this meta-reasoning is that:

(*xiii*)   ~PA proves: [(A*y*)(*Q*(*u*, *y*) => (E*z*)(*z*=<*y* & *S*(*u*, *z*)))]

   ... *Assuming PA is consistent*

(*b*) The second half of the argument is essentially:

(*i*)   PA proves: [~(A*y*)(*Q*(*u*, *y*) => (E*z*)(*z*=<*y* & *S*(*u*, *z*)))]

   ... *Hypothesis*

(*ii*)   PA proves: [~(A*y*)(*Q*(*u*, *y*) => (E*z*)(*z*=<*y* & *S*(*u*, *z*)))]

   => (E*r*)PA proves: [*s*(*u*, *r*)]

   ... *By the Self-reference meta-Lemma §7.2(v)*

(*iii*)   (E*r*)PA proves: [*s*(*u*, *r*)]

   ... *From (i) and (ii) by Modus Ponens in the meta-theory*

(*iv*)   (A*n*)PA proves: [*s*(*u*, *n*) <=> *S*(*u*, *n*)]

   ... *By the Representation meta-Lemmas §7.2(i)*

(*v*)   (E*r*)PA proves: [*S*(*u*, *r*)]

   ... *From (iii) and (iv) by Modus Ponens, assuming (i)*

(*vi*)    (A*n*)PA proves: [~*q*(*u*, *n*)]

   ... *From (i) and §7.2(iii) assuming PA is consistent*

(*vii*)   (A*n*)PA proves: [~*q*(*u*, *n*) <=> ~*Q*(*u*, *n*)]

   ... *By the Representation meta-Lemmas §4.1(i)*

(*viii*)  (A*n*)PA proves: [~*Q*(*u*, *n*)]

   ... *From (vi) and (vii) by Modus Ponens, assuming (i)*

(*ix*)    (A*n*)PA proves: [(A*z*)(*z*=<*n* => ~*Q*(*u*, *z*))]

   ... *From (viii) by the logical axioms of PA, assuming (i)*

   (*x*)(*1*)   (E*r*)PA proves: [*r*=<*y*]

      ... *Theorem of PA*

   (*x*)(*2*)   (E*r*)PA proves: [*S*(*u*, *r*)]

      ... *From (v) above, assuming (i)*

   (*x*)(*3*)   (E*r*)PA proves: [*r*=<*y* & *S*(*u*, *r*)]

      ... *From (1) and (2), assuming (i) – Invalid deduction*

   (*x*)(*4*)   PA proves: [(E*z*)(*z*=<*y* & *S*(*u*, *z*))]

      ... *From (3), assuming (i)*

(*xi*)    (E*r*)PA proves: [*r*=<*y* => (E*z*)(*z*=<*y* & *S*(*u*, *z*))]

   ... *From (x)(1) and (x)(4) by the Deduction Theorem, assuming (i)*

(*xii*)   (A*r*)PA proves: [*y*=<*r* v *r*=<*y*]

    *... Theorem of PA*

(*xiii*)   PA proves: [~*Q*(*u*, *y*) v (E*z*)(*z*=<*y* & *S*(*u*, *z*))]

    *... From (ix), (xi) and (xii), by the logical axioms of PA*

(*xiv*)   PA proves: [(A*y*)(*Q*(*u*, *y*) => (E*z*)(*z*=<*y* & *S*(*u*, *z*)))]

    *... From (xiii), by the logical axioms of PA, Modus Ponens and Generalisation, assuming (i)*

Since (*xiv*) contradicts (*i*), the classical conclusion drawn from the above is that:

(*xv*)   ~PA proves: [~(A*y*)(*Q*(*u*, *y*) => (E*z*)(*z*=<*y* & *S*(*u*, *z*)))]

    *... Assuming that PA is consistent*

Now we note that the above meta-reasoning is based on the invalid deduction (*x*)(3). This appeals to the implicitly Platonistic assumption that, if we assume that PA proves: [~(A*y*)(*Q*(*u*, *y*) => (E*z*)(*z*=<*y* & *S*(*u*, *z*)))], then we can also conclude that PA proves: [*S*(*u*, *r*)] for some unspecified, but specific, natural number *r*, as in §7.3(*b*)(*i*)(*2*).

However, from the assertion "PA proves: [~(A*y*)(*Q*(*u*, *y*) => (E*z*)(*z*=<*y* & *S*(*u*, *z*)))]" and Rosser's Self-reference meta-Lemmas §7.2(*ii*)-(*iv*), we may constructively only conclude the non-specific meta-mathematical assertion "(E*r*)PA proves: [*S*(*u*, *r*)]", as in (*v*) above. We may not further introduce "PA proves: [*S*(*u*, *r*)]" for a specific *r*, as in §7.3(*b*)(*i*)(*2*), except as an additional premise[12] in the application of the Deduction Theorem[13].

---

[12] We note, however, that for a given *r*, "PA proves: *S*(*u*, *r*)" can not be treated as a meaningful premise, since it is meta-mathematically PA-decidable.

However, this would not then yield the consequence (*xi*) as Rosser intended. So the formal deduction of §7.3(*b*)(*ii*) from §7.3(*b*)(*i*)(1)-(4) can reasonably be held as invalid in any non-semantic meta-reasoning such as (*x*)(1)-(4) above.

## 8. Conclusion

(*a*) We can thus reasonably conclude that the classical interpretations of Gödel's and Rosser's meta-arguments do not unequivocally establish undecidable sentences in any formal system in which:

  (*i*)   The proper *Axioms* and *Rules of Inference* are recursively definable.

  (*ii*)  Every recursive relation is representable.

(*b*) We can also reasonably conclude further that:

  (*i*)   Gödel's undecidable sentence GUS is actually decidable in PA under a reasonable interpretation of his meta-reasoning that, avoiding any appeal to

---

[13] Rosser's meta-argument seems to be that, assuming [$H(\underline{h})$] is provable in PA, and that $j$ is the Gödel-number of a proof of [$H(\underline{h})$] in PA, we can express this meta-mathematically as:

  (*i*)   PA proves: [$H(\underline{h})$] => PA proves: [$q(\underline{h}, \underline{j})$].

From this he apparently concludes by the Deduction Theorem that:

  (*ii*)  PA proves: [$H(\underline{h}) => q(\underline{h}, \underline{j})$].

However, (*ii*) is invalid under interpretation in IA. Thus, given any natural number $h$, we could construct a Turing-machine that would decompose $h$ to check whether it is in fact the Goedel-number of some well-formed formula [$H(x)$] of PA. Then, assuming the provability of [$H(\underline{h})$], the program would conclude that the meta-assertion (*ii*) is invalid as it does not hold for every natural number $j$.

I thus argue that the formal meta-mathematical expression of (*i*) is the meta-statement:

  (*iii*) PA proves: [$H(\underline{h})$] => (E$j$)PA proves: [$q(\underline{h}, \underline{j})$].

However, as I argue further, the Deduction Theorem cannot be applied to the formal meta-assertion (*iii*) so as to yield Rosser's undecidable proposition RUS in PA.

30the truth of propositions of IA, is constructive and intuitionistically unobjectionable;

(*ii*)   PA is not omega-consistent under such interpretation;

(*iii*)  If we take IA to be the standard model for PA, this implies that standard PA is semantically inconsistent under the standard interpretations of Gödel's meta-reasoning;

(*iv*)   Rosser's extension of Gödel's meta-reasoning, establishing an undecidable proposition RUS in a consistent PA, is invalid under such an interpretation.

## References

[An01]   Anand, Bhupinder Singh. 2001. *Beyond Gödel: Simply constructive systems of first order Peano's Arithmetic that do not yield undecidable propositions by Gödel's reasoning*. Alix Comsi, Mumbai (*Unpublished*).

<*Web page*: http://alixcomsi.com/Constructivity_abstracts.htm>

[An02]   Anand, Bhupinder Singh. 2002. *Paradox regained: Life beyond Gödel's shadow*. Alix Comsi, Mumbai (*Unpublished*).

<*Web page*: http://alixcomsi.com/Constructivity_preamble_Rev1.htm>

[Go31]   Gödel, Kurt. 1931. *On formally undecidable propositions of Principia Mathematica and related systems I*. Also in M. Davis (ed.). 1965. The Undecidable. Raven Press, New York.

<*Web version*: http://www.ddc.net/ygg/etext/godel/godel3.htm>

[Me64]   Mendelson, Elliott. 1964. Introduction to Mathematical Logic. Van Norstrand, Princeton.

<*Author's home page*: http://sard.math.qc.edu/Web/Faculty/mendelso.htm>(*Acknowledgement: I am grateful to Catherine Christer Hennix, Damjan Bojadziev, Andrew Boucher, Dennis E. Hamilton, Charlie Volkstorf, Aatu Koskensilta and other correspondents for their constructive comments on various issues addressed in this paper.*)

(*Updated: Friday 9$^{th}$ May 2003 8:56:40 AM by re@alixcomsi.com*)